# CONVEXITY, TRANSLATION INVARIANCE AND SUBADDITIVITY FOR $G$-EXPECTATIONS AND RELATED RISK MEASURES

By Long Jiang[1]

*China University of Mining and Technology, Fudan University*

Under the continuous assumption on the generator $g$, Briand et al. [*Electron. Comm. Probab.* **5** (2000) 101–117] showed some connections between $g$ and the conditional $g$-expectation $(\mathcal{E}_g[\cdot|\mathcal{F}_t])_{t\in[0,T]}$ and Rosazza Gianin [*Insurance: Math. Econ.* **39** (2006) 19–34] showed some connections between $g$ and the corresponding dynamic risk measure $(\rho_t^g)_{t\in[0,T]}$. In this paper we prove that, without the additional continuous assumption on $g$, a $g$-expectation $\mathcal{E}_g$ satisfies translation invariance if and only if $g$ is independent of $y$, and $\mathcal{E}_g$ satisfies convexity (resp. subadditivity) if and only if $g$ is independent of $y$ and $g$ is convex (resp. subadditive) with respect to $z$. By these conclusions we deduce that the static risk measure $\rho^g$ induced by a $g$-expectation $\mathcal{E}_g$ is a convex (resp. coherent) risk measure if and only if $g$ is independent of $y$ and $g$ is convex (resp. sublinear) with respect to $z$. Our results extend the results in Briand et al. [*Electron. Comm. Probab.* **5** (2000) 101–117] and Rosazza Gianin [*Insurance: Math. Econ.* **39** (2006) 19–34] on these subjects.

## 1. Introduction and preliminaries.

1.1. *Introduction.* By Pardoux and Peng [16] we know that there exists a unique adapted and square integrable solution to a backward stochastic differential equation (BSDE in short) of type

$$(1.1) \qquad y_t = \xi + \int_t^T g(s, y_s, z_s)\, ds - \int_t^T z_s \cdot dB_s, \qquad 0 \le t \le T,$$

Received September 2005; revised June 2007.

[1]Supported by NSFC Grant 10671205 and CPSF Grant 20060400158 and National Basic Research Program of China Grant 2007CB814901.

*AMS 2000 subject classifications.* Primary 60H10; secondary 60H30, 91B30.

*Key words and phrases.* Backward stochastic differential equation, $g$-expectation, translation invariance of $g$-expectation, convexity of $g$-expectation, risk measure.

providing that the function $g$ is Lipschitz in both variables $y$ and $z$, and $\xi$ and $(g(t,0,0))_{t\in[0,T]}$ are square integrable. $g$ is called the generator of the BSDE (1.1) and $(g,T,\xi)$ are called the parameters of (1.1). We denote the unique solution of (1.1) by $(Y_t(g,T,\xi), Z_t(g,T,\xi))_{t\in[0,T]}$. When $g$ also satisfies $g(\cdot,y,0)\equiv 0$ for any $y$, then $Y_0(g,T,\xi)$, denoted by $\mathcal{E}_g[\xi]$, is called the $g$-expectation of $\xi$; $Y_t(g,T,\xi)$, denoted by $\mathcal{E}_g[\xi|\mathcal{F}_t]$, is called the conditional $g$-expectation of $\xi$; see Peng [17].

$g$-expectation is a kind of nonlinear expectation. The original motivation for studying $g$-expectation comes from the theory of expected utility, which is the foundation of modern mathematical economics and is challenged by the well-known Allais paradox. Since the notion of $g$-expectation was introduced, many properties of $g$-expectation have been studied in [3, 4, 5, 6, 13, 14, 15, 17, 18, 19]. Chen and Epstein [5] gave an application of $g$-expectation to recursive utility. Coquet et al. [6] obtained a very interesting result. They proved that if a filtration consistent (nonlinear) expectation $\mathcal{E}$ can be dominated by a kind of $g$-expectation, then $\mathcal{E}$ must be a $g$-expectation. More recently, Rosazza Gianin [20, 21] first introduced some examples of risk measures via $g$-expectations and conditional $g$-expectations:

$$(1.2) \qquad \rho^g(\xi) := \mathcal{E}_g[-\xi], \qquad \rho_t^g(\xi) := \mathcal{E}_g[-\xi|\mathcal{F}_t] \qquad \forall t\in[0,T].$$

Under an additional continuity assumption (A4) (see Section 1.2), with the help of Proposition 2.3 of Briand et al. [3], Rosazza Gianin [21] showed us that $(\rho_t^g)_{t\in[0,T]}$ is a dynamic convex (resp. coherent) risk measure if and only if $g$ is independent of $y$ and is convex (resp. sublinear) with respect to $z$. Barrieu and El Karoui [2] and Peng [19] also obtained some results on this subject.

The main objective of this paper is to explore some fundamental characteristics of $g$-expectations which are related to risk measures. The major contributions of this paper are:

(a) We establish a general Representation Lemma for generators of BSDEs under the usual assumptions (A1) and (A2), which generalizes Proposition 2.3 of [3] and helps us to confirm the same necessary and sufficient conditions in [3] and [21] without the additional continuity assumption (A4). We hope that it turns out to be useful in other situations, as well.

(b) Under the usual assumptions (A1) and (A3), without any additional assumptions on $g$, we prove that if the static risk measure $\rho^g$, which is an operator, is a convex (resp. coherent) risk measure, then the corresponding dynamic risk measure $(\rho_t^g)_{t\in[0,T]}$, which is an operator system, is a dynamic convex (resp. coherent) risk measure, and the generator $g$ is independent of $y$ and is convex (resp. sublinear) with respect to $z$.

The remainder of this paper is organized as follows. In Section 1.2, we introduce some preliminaries. In Section 2, we establish a general Representation Lemma for generators of BSDEs. In Section 3, under the usual assumptions (A1) and (A3), we obtain some necessary and sufficient conditions for translation invariance, convexity, subadditivity and positive homogeneity of $g$-expectations, respectively. In Section 4, we state our results on static risk measure $\rho^g$ and dynamic risk measure $(\rho_t^g)_{t \in [0,T]}$.

1.2. *Preliminaries.* Let $T > 0$ be a fixed time horizon; let $(\Omega, \mathcal{F}, P)$ be a probability space and $(B_t)_{t \geq 0}$ be a $d$-dimensional standard Brownian motion on this space such that $B_0 = 0$; let $(\mathcal{F}_t)_{t \geq 0}$ be the augmented natural filtration generated by $(B_t)_{t \geq 0}$ and satisfy the usual conditions. Let $\mathcal{M}_{\mathcal{F}}(\mathbf{R}^n)$ denote the space of all $\mathbf{R}^n$-valued, $(\mathcal{F}_t)$-progressively measurable processes. We set $\mathcal{H}_{\mathcal{F}}^2(0,T;\mathbf{R}^n) := \{\psi \in \mathcal{M}_{\mathcal{F}}(\mathbf{R}^n); \|\psi\|_2^2 := \mathbf{E}[\int_0^T |\psi_t|^2 \, dt] < \infty\}$, $\mathcal{S}_{\mathcal{F}}^2(0,T;\mathbf{R}) := \{\psi \in \mathcal{M}_{\mathcal{F}}(\mathbf{R}); \psi$ is continuous and $\mathbf{E}[\sup_{0 \leq t \leq T} |\psi_t|^2] < \infty\}$, $L^2(\mathcal{F}_t) := \{\xi; \xi \text{ is } \mathbf{R}\text{-valued}, \mathcal{F}_t\text{-measurable random variable}, \mathbf{E}[\xi^2] < \infty\}$.

The generator $g$ of a BSDE is a function $g : \Omega \times [0,T] \times \mathbf{R} \times \mathbf{R}^d \longmapsto \mathbf{R}$ such that $(g(t,y,z))_{t \in [0,T]}$ is progressively measurable for each $(y,z) \in \mathbf{R} \times \mathbf{R}^d$, and $g$ also satisfies the following usual assumptions (A1) and (A2):

(A1) There exists a constant $K \geq 0$ such that $dP \times dt$-a.s., $\forall y_1, y_2, z_1, z_2$, $|g(t,y_1,z_1) - g(t,y_2,z_2)| \leq K(|y_1 - y_2| + |z_1 - z_2|)$.

(A2) The process $(g(t,0,0))_{t \in [0,T]} \in \mathcal{H}_{\mathcal{F}}^2(0,T;\mathbf{R})$.

(A3) $dP \times dt$-a.s., for any $y \in \mathbf{R}$, $g(\cdot, y, 0) \equiv 0$.

(A4) $P$-a.s., for any $y \in \mathbf{R}, z \in \mathbf{R}^d, t \mapsto g(t,y,z)$ is continuous.

Let (A1) and (A2) hold for $g$. By [16], for each $\xi \in L^2(\mathcal{F}_T)$, (1.1) has a unique solution in $\mathcal{S}_{\mathcal{F}}^2(0,T;\mathbf{R}) \times \mathcal{H}_{\mathcal{F}}^2(0,T;\mathbf{R}^d)$, which is denoted by

$$(Y_t(g,T,\xi), Z_t(g,T,\xi))_{t \in [0,T]}.$$

We recall the notions of $g$-expectation and conditional $g$-expectation and some properties given in Peng [17]. In the following Definitions 1.1 and 1.2 and Lemma 1.1, $g$ is assumed to satisfy (A1) and (A3).

DEFINITION 1.1 ([17]). The $g$-expectation $\mathcal{E}_g[\cdot] : L^2(\mathcal{F}_T) \longmapsto \mathbf{R}$ is defined by $\mathcal{E}_g[\xi] := Y_0(g,T,\xi)$.

DEFINITION 1.2 ([17]). The conditional $g$-expectation of $\xi$ with respect to $\mathcal{F}_t$ is defined by $\mathcal{E}_g[\xi|\mathcal{F}_t] := Y_t(g,T,\xi)$, which is the unique random variable $\eta$ in $L^2(\mathcal{F}_t)$ such that $\mathcal{E}_g[\xi \mathbf{1}_A] = \mathcal{E}_g[\eta \mathbf{1}_A]$, for all $A \in \mathcal{F}_t$.

LEMMA 1.1 ([17]). (i) *For each constant* $c, \mathcal{E}_g[c] = c$.

(ii) *If* $X_1 \geq X_2, a.s.,$ *then* $\mathcal{E}_g[X_1] \geq \mathcal{E}_g[X_2]$.

(iii) *If* $X_1 \geq X_2, a.s.,$ *and* $P(X_1 > X_2) > 0,$ *then* $\mathcal{E}_g[X_1] > \mathcal{E}_g[X_2]$.

(iv) *If* $g$ *is independent of* $y$, *that is,* $g$ *is defined on* $\Omega \times [0,T] \times \mathbf{R}^d$, *then* $\mathcal{E}_g[X + \eta|\mathcal{F}_t] = \mathcal{E}_g[X|\mathcal{F}_t] + \eta, \forall X \in L^2(\mathcal{F}_T), \eta \in L^2(\mathcal{F}_t)$.

**2. Representation lemma for generators of BSDEs.** For studying a kind of converse comparison problem, Proposition 2.3 in Briand et al. [3] showed us that for any $(y, z) \in \mathbf{R} \times \mathbf{R}^d$ and $t \in [0, T[$, the equality

$$g(t, y, z) = L^2 - \lim_{\varepsilon \to 0^+} \frac{1}{\varepsilon} [Y_t(g, t + \varepsilon, y + z \cdot (B_{t+\varepsilon} - B_t)) - y]$$

holds under (A1), (A2), (A4) and $\mathbf{E}[\sup_{0 \le t \le T} |g(t, 0, 0)|^2] < \infty$. For studying Jensen's inequality for $g$-expectation, [15] got the following Proposition 2.1.

PROPOSITION 2.1 (Theorem 3.3 in [15]). *Let* (A1) *and* (A2) *hold for* $g$; *let* $1 \le p \le 2$. *Then for any triplet* $(t, y, z) \in [0, T[ \times \mathbf{R} \times \mathbf{R}^d$, *the following two statements are equivalent:*

(i) $g(t, y, z) = L^p - \lim_{\varepsilon \to 0^+} \frac{1}{\varepsilon} [Y_t(g, t + \varepsilon, y + z \cdot (B_{t+\varepsilon} - B_t)) - y]$.

(ii) $g(t, y, z) = L^p - \lim_{\varepsilon \to 0^+} \mathbf{E}[\frac{1}{\varepsilon} \int_t^{t+\varepsilon} g(s, y, z) \, ds | \mathcal{F}_t]$.

Further studying shows that many problems on BSDEs are related to this kind of representation problem. In this section, we will establish a general Representation Lemma for generators of BSDEs under (A1) and (A2), which generalizes Proposition 2.3 of [3] and will be used frequently.

LEMMA 2.1 (Representation lemma). *Let* (A1) *and* (A2) *hold for* $g$. *Let* $1 \le p < 2$. *Then for each* $(y, z) \in \mathbf{R} \times \mathbf{R}^d$, *the equality*

$$g(t, y, z) = L^p - \lim_{\varepsilon \to 0^+} \frac{1}{\varepsilon} [Y_t(g, t + \varepsilon, y + z \cdot (B_{t+\varepsilon} - B_t)) - y]$$

*holds for almost every* $t \in [0, T[$.

In order to prove Lemma 2.1, we introduce the following proposition.

PROPOSITION 2.2. *Let* $q > 1$; *let* $1 \le p < q$. *Set* $\mathcal{H}_{\mathcal{F}}^q(0, T; \mathbf{R}) := \{\psi \in \mathcal{M}_F(\mathbf{R}); \mathbf{E}[\int_0^T |\psi_t|^q \, dt] < \infty\}$. *Then for any* $\psi \in \mathcal{H}_{\mathcal{F}}^q(0, T; \mathbf{R})$, *we have*

$$\psi_t = L^p - \lim_{\varepsilon \to 0^+} \frac{1}{\varepsilon} \int_t^{t+\varepsilon} \psi_s \, ds \qquad a.e. \ t \in [0, T[.$$

PROOF. Since $\psi \in \mathcal{H}_{\mathcal{F}}^q(0, T; \mathbf{R})$, the Fubini theorem yields $\int_0^T \mathbf{E}[|\psi_t|^q] \, dt = \mathbf{E}[\int_0^T |\psi_t|^q \, dt] < \infty$. Thus $\mathbf{E}[|\psi_t|^q] < \infty$, a.e. $t \in [0, T]$. By the Lebesgue lemma (see Lemma 18.4 of [12]), we know that the equality

$$\lim_{\varepsilon \to 0^+} \frac{1}{\varepsilon} \int_t^{t+\varepsilon} \mathbf{E}[|\psi_s|^q] \, ds = \mathbf{E}[|\psi_t|^q]$$

holds for almost every $t \in [0, T[$.

Also by $\psi \in \mathcal{H}_{\mathcal{F}}^q(0, T; \mathbf{R})$ we understand that $|\int_0^T \psi_t \, dt| < \infty$, a.s. Therefore, by the Lebesgue lemma we have

$$(2.1) \qquad \lim_{\varepsilon \to 0^+} \frac{1}{\varepsilon} \int_t^{t+\varepsilon} \psi_s \, ds = \psi_t \qquad \text{a.e., a.s.}$$

Hence

$$\lim_{\varepsilon \to 0^+} \frac{1}{\varepsilon} \int_t^{t+\varepsilon} \psi_s \, ds = \psi_t \qquad \text{a.s., a.e.}$$

Thus there exists a subset $S \subseteq [0, T[$ such that the Lebesgue measure $\lambda([0, T] \setminus S)$ of $[0, T] \setminus S$ equals 0, and for each $t \in S$ we have

$$(2.2) \qquad \lim_{\varepsilon \to 0^+} \frac{1}{\varepsilon} \int_t^{t+\varepsilon} \psi_s \, ds = \psi_t \qquad \text{a.s.,}$$

$$(2.3) \qquad \mathbf{E}[|\psi_t|^q] < \infty, \qquad \lim_{\varepsilon \to 0^+} \frac{1}{\varepsilon} \int_t^{t+\varepsilon} \mathbf{E}[|\psi_s|^q] \, ds = \mathbf{E}[|\psi_t|^q].$$

For any $t \in S$, by (2.3) we know that there exists a constant $\delta_t > 0$ such that

$$(2.4) \qquad \frac{1}{\varepsilon} \int_t^{t+\varepsilon} \mathbf{E}[|\psi_s|^q] \, ds \leq \mathbf{E}[|\psi_t|^q] + 1 \qquad \forall \varepsilon \in \, ]0, \delta_t].$$

For any $t \in S$, $\varepsilon \in \, ]0, \delta_t]$, we set $X_t^\varepsilon := |\frac{1}{\varepsilon} \int_t^{t+\varepsilon} \psi_s \, ds|$. Then for any $N > 0$, by Hölder's inequality, Fubini's theorem and (2.4) we have

$$\int_{\{X_t^\varepsilon > N\}} \left| \frac{1}{\varepsilon} \int_t^{t+\varepsilon} \psi_s \, ds \right|^p dP \leq \int_{\{X_t^\varepsilon > N\}} \frac{1}{N^{q-p}} \left| \frac{1}{\varepsilon} \int_t^{t+\varepsilon} \psi_s \, ds \right|^q dP$$

$$\leq \int_{\{X_t^\varepsilon > N\}} \frac{1}{N^{q-p}} \left[ \frac{1}{\varepsilon} \int_t^{t+\varepsilon} |\psi_s|^q \, ds \right] dP$$

$$\leq \frac{1}{N^{q-p}} \mathbf{E} \left[ \frac{1}{\varepsilon} \int_t^{t+\varepsilon} |\psi_s|^q \, ds \right]$$

$$\leq \frac{1}{N^{q-p}} [\mathbf{E}[|\psi_t|^q] + 1].$$

Thus $\{|\frac{1}{\varepsilon} \int_t^{t+\varepsilon} \psi_s \, ds|^p; \varepsilon \in \, ]0, \delta_t]\}$ are uniformly integrable. Combining this conclusion with (2.2), we conclude that for each $t \in S$, we have

$$(2.5) \qquad \psi_t = L^p - \lim_{\varepsilon \to 0^+} \frac{1}{\varepsilon} \int_t^{t+\varepsilon} \psi_s \, ds.$$

The proof of Proposition 2.2 is complete. $\square$

PROOF OF LEMMA 2.1. Since $(g(t, 0, 0))_{t \in [0, T]} \in \mathcal{H}_{\mathcal{F}}^2(0, T; \mathbf{R})$ and $g$ satisfies (A1), we know that for each $(y, z) \in \mathbf{R} \times \mathbf{R}^d$, $(g(t, y, z))_{t \in [0, T]} \in \mathcal{H}_{\mathcal{F}}^2(0, T; \mathbf{R})$.

Then for any $1 \leq p < 2$ and any $(y, z) \in \mathbf{R} \times \mathbf{R}^d$, Proposition 2.2 and Jensen's inequality yield

$$(2.6) \quad g(t, y, z) = L^p - \lim_{\varepsilon \to 0^+} \mathbf{E} \left[ \frac{1}{\varepsilon} \int_t^{t+\varepsilon} g(s, y, z) \, ds \Big| \mathcal{F}_t \right] \qquad \text{a.e. } t \in [0, T[.$$

Thus Lemma 2.1 follows from (2.6) and Proposition 2.1 immediately. $\quad \square$

REMARK 2.1. Consider a financial market where derivatives pricing is constrained by BSDEs with generator $g$; Lemma 2.1 may help us to find the pricing mechanism, that is, the function $g$.

**3. Translation invariance, convexity, subadditivity and positive homogeneity for $g$-expectations.** In this section, we study some properties of $g$-expectations such as translation invariance, convexity, subadditivity and positive homogeneity. All these properties are related to risk measures via $g$-expectations. We obtain some necessary and sufficient conditions on these problems, respectively. The main differences between our results and other known results on these problems such as those results in Briand et al. [3] and Rosazza Gianin [21] are:

(a) We can use the $g$-expectation $\mathcal{E}_g[\cdot]$, which is an operator, to describe the character of the generator $g$; on the other hand, [3] always used the conditional $g$-expectation $(\mathcal{E}_g[\cdot|\mathcal{F}_t])_{t \in [0,T]}$, which is an operator system, to describe the character of $g$.

(b) Our results are obtained under the usual assumptions (A1) and (A3); on the other hand, the necessary and sufficient conditions given in [3] and [21] are always obtained under the assumptions (A1), (A3) and the additional continuity assumption (A4).

From now on, for any pair $(y, z) \in \mathbf{R} \times \mathbf{R}^d$, we set

$$S_y^z(g) := \Big\{ t \in [0, T[ \Big|$$

$$g(t, y, z) = L^1 - \lim_{\varepsilon \to 0^+} \frac{1}{\varepsilon} [Y_t(g, t+\varepsilon, y + z \cdot (B_{t+\varepsilon} - B_t)) - y] \Big\}.$$

If $g$ is independent of $y$, then for any $z \in \mathbf{R}^d$, we set

$$S^z(g) := \Big\{ t \in [0, T[ \Big| g(t, z) = L^1 - \lim_{\varepsilon \to 0^+} \frac{1}{\varepsilon} Y_t(g, t+\varepsilon, z \cdot (B_{t+\varepsilon} - B_t)) \Big\}.$$

3.1. *Translation invariance for $g$-expectation.* If $g$ is independent of $y$, then by Lemma 1.1(iv) we know that the $g$-expectation $\mathcal{E}_g$ satisfies translation invariance. We now investigate the inverse problem. We have the following theorem.

THEOREM 3.1 (Translation invariance for *g*-expectation). *Let* (A1) *and* (A3) *hold for g. Then the following three statements are equivalent:*

(i) $\mathcal{E}_g[\xi + c] = \mathcal{E}_g[\xi] + c, \forall \xi \in L^2(\mathcal{F}_T),\ c \in \mathbf{R}.$   (Translation invariance.)

(ii) *For* $\forall \xi \in L^2(\mathcal{F}_T), c \in \mathbf{R}, \forall t \in [0, T],$

$$\mathcal{E}_g[\xi + c | \mathcal{F}_t] = \mathcal{E}_g[\xi | \mathcal{F}_t] + c, \qquad P\text{-}a.s.$$

(iii) *g is independent of y.*

PROOF. (iii) $\Rightarrow$ (ii) follows from Lemma 1.1. (ii) $\Rightarrow$ (i) is trivial. Now let us prove that (i) $\Rightarrow$ (iii). Suppose that (i) holds.

For any $c \in \mathbf{R}$, we define a new generator

$$g^c(t, y, z) := g(t, y - c, z) \qquad \forall t \in [0, T],\ y \in \mathbf{R},\ z \in \mathbf{R}^d.$$

Then $g^c$ satisfies (A1), (A2) and (A3).

For any $X \in L^2(\mathcal{F}_T)$, by the uniqueness of solution of BSDE we can verify easily that

$$(Y_t(g^c, T, X + c), Z_t(g^c, T, X + c))_{t \in [0,T]} = (Y_t(g, T, X) + c, Z_t(g, T, X))_{t \in [0,T]}.$$

It follows that

$$\mathcal{E}_{g^c}[X + c] = Y_0(g^c, T, X + c) = Y_0(g, T, X) + c = \mathcal{E}_g[X] + c.$$

Combining the above equality with (i) we have

$$\mathcal{E}_{g^c}[X + c] = \mathcal{E}_g[X + c] \qquad \forall X \in L^2(\mathcal{F}_T).$$

Hence for any given $c \in \mathbf{R}$, we have

$$(3.1) \qquad \mathcal{E}_{g^c}[\xi] = \mathcal{E}_g[\xi] \qquad \forall \xi \in L^2(\mathcal{F}_T).$$

It follows from (3.1) and Proposition 3.4 of [13] that for any $\xi \in L^2(\mathcal{F}_T)$, we have

$$(3.2) \qquad P\text{-a.s.}, \forall t \in [0, T] \qquad \mathcal{E}_{g^c}[\xi | \mathcal{F}_t] = \mathcal{E}_g[\xi | \mathcal{F}_t].$$

Then for any $(y, z) \in \mathbf{R} \times \mathbf{R}^d$ and for any $t \in S_y^z(g^c) \cap S_y^z(g)$, (3.2) yields

$$(3.3) \qquad P\text{-a.s.},\ g^c(t, y, z) = g(t, y, z).$$

By the representation lemma we understand that

$$(3.4) \qquad \lambda([0, T] \setminus (S_y^z(g^c) \cap S_y^z(g))) = 0,$$

where $\lambda$ denotes the Lebesgue measure. It follows from (3.3) and (3.4) that

$$(3.5) \qquad dP \times dt\text{-a.s.}, \qquad g^c(t, y, z) = g(t, y, z).$$

Since $g$ and $g^c$ are both Lipschitz with respect to $(y, z)$, it follows that

$$(3.6) \qquad dP \times dt\text{-a.s.},\ \forall y \in \mathbf{R}, z \in \mathbf{R}^d \qquad g^c(t, y, z) = g(t, y, z);$$

that is, for any given $c \in \mathbf{R}$, we have $g^c = g$. Thus for any $y \in \mathbf{R}$, we have

$$(3.7) \qquad dP \times dt\text{-a.s.}, \ \forall z \in \mathbf{R}^d \qquad g(t, y, z) = g(t, 0, z).$$

Therefore (iii) follows from (3.7) and the Lipschitz assumption (A1). $\square$

REMARK 3.1. Under (A1), (A3) and (A4), Briand et al. [3] proved that (ii) is equivalent to (iii) in Theorem 3.1; see Lemmas 4.2 and 4.3 in [3].

3.2. *Convexity, subadditivity and positive homogeneity for g-expectations.* For studying a control problem, El Karoui, Peng and Quenez [8] studied concave BSDEs. For studying dynamic risk measures, Rosazza Gianin [20, 21] studied the convexity, subadditivity and positive homogeneity of conditional $g$-expectations. The reader can see some results of [21] in Remark 3.2. Now let us introduce our results.

THEOREM 3.2 (Convexity for g-expectation). *Let* (A1) *and* (A3) *hold for* $g$. *Then the following three statements are equivalent:*

(i) $\mathcal{E}_g[\cdot]$ *is convex.* (Convexity.)

(ii) *For any* $t \in [0, T]$, $\mathcal{E}_g[\cdot | \mathcal{F}_t]$ *is convex, that is,* $\forall \xi, \eta \in L^2(\mathcal{F}_T), \alpha \in [0, 1]$,

$$\mathcal{E}_g[\alpha\xi + (1-\alpha)\eta | \mathcal{F}_t] \leq \alpha\mathcal{E}_g[\xi | \mathcal{F}_t] + (1-\alpha)\mathcal{E}_g[\eta | \mathcal{F}_t], \qquad P\text{-a.s.}$$

(iii) $g$ *is independent of* $y$ *and* $g$ *is convex with respect to* $z$, *that is, for any* $z_1, z_2 \in \mathbf{R}^d, \alpha \in [0, 1]$,

$$g(t, \alpha z_1 + (1-\alpha)z_2) \leq \alpha g(t, z_1) + (1-\alpha)g(t, z_2), \qquad dP \times dt\text{-a.s.}$$

PROOF. (iii) $\Rightarrow$ (ii) follows from the well-known comparison theorem; the argument is analogous to the argument of Proposition 3.5 in El Karoui, Peng and Quenez [8] when those authors studied concave BSDEs. (ii) $\Rightarrow$ (i) is trivial. Now let us prove that (i) $\Rightarrow$ (iii).

Suppose that (i) holds. First, let us prove that the convexity of $g$-expectation implies the translation invariance. Indeed, for any $\xi \in L^2(\mathcal{F}_T), c \in \mathbf{R}, \alpha \in [0, 1]$, by (i) and Lemma 1.1 we have

$$(3.8) \quad \mathcal{E}_g[\alpha\xi + (1-\alpha)c] \leq \alpha\mathcal{E}_g[\xi] + (1-\alpha)\mathcal{E}_g[c] = \alpha\mathcal{E}_g[\xi] + (1-\alpha)c.$$

Thus for any $\xi \in L^2(\mathcal{F}_T), c \in \mathbf{R}$ and any positive integer $n$, we have

$$\mathcal{E}_g\left[\left(1 - \frac{1}{n}\right)\xi + c\right] = \mathcal{E}_g\left[\left(1 - \frac{1}{n}\right)\xi + \frac{1}{n}(nc)\right]$$

$$\leq \left(1 - \frac{1}{n}\right)\mathcal{E}_g[\xi] + c.$$

Since the operator $\mathcal{E}_g[\cdot]$ is continuous in $L^2$ sense, we have

$$\mathcal{E}_g[\xi + c] = \lim_{n \to \infty} \mathcal{E}_g\left[\left(1 - \frac{1}{n}\right)\xi + c\right]$$

$$\leq \lim_{n \to \infty} \left(1 - \frac{1}{n}\right)\mathcal{E}_g[\xi] + c$$

$$= \mathcal{E}_g[\xi] + c.$$

Hence we have

$$(3.9) \qquad \mathcal{E}_g[\xi + c] \leq \mathcal{E}_g[\xi] + c \qquad \forall \xi \in L^2(\mathcal{F}_T), \ c \in \mathbf{R}.$$

Therefore

$$(3.10) \qquad \mathcal{E}_g[\xi] = \mathcal{E}_g[\xi + c - c] \leq \mathcal{E}_g[\xi + c] - c \qquad \forall \xi \in L^2(\mathcal{F}_T), \ c \in \mathbf{R}.$$

It follows from the above two inequalities that

$$(3.11) \qquad \mathcal{E}_g[\xi + c] = \mathcal{E}_g[\xi] + c \qquad \forall \xi \in L^2(\mathcal{F}_T), \ c \in \mathbf{R}.$$

Thus the $g$-expectation $\mathcal{E}_g$ satisfies the translation invariance. Then by Theorem 3.1 we conclude that $g$ is independent of $y$.

Second, let us prove that for each $\xi, \eta \in L^2(\mathcal{F}_T), \alpha \in [0, 1], P$-a.s.,

$$(3.12) \qquad \mathcal{E}_g[\alpha\xi + (1-\alpha)\eta | \mathcal{F}_t] \leq \alpha\mathcal{E}_g[\xi | \mathcal{F}_t] + (1-\alpha)\mathcal{E}_g[\eta | \mathcal{F}_t] \qquad \forall t \in [0, T].$$

We set

$$A := \{\mathcal{E}_g[\alpha\xi + (1-\alpha)\eta | \mathcal{F}_t] > \alpha\mathcal{E}_g[\xi | \mathcal{F}_t] + (1-\alpha)\mathcal{E}_g[\eta | \mathcal{F}_t]\}.$$

Then $A \in \mathcal{F}_t$. Suppose by contradiction that $P(A) > 0$. Then

$$\mathbf{1}_A\mathcal{E}_g[\alpha\xi + (1-\alpha)\eta | \mathcal{F}_t] - \mathbf{1}_A(\alpha\mathcal{E}_g[\xi | \mathcal{F}_t] + (1-\alpha)\mathcal{E}_g[\eta | \mathcal{F}_t]) \geq 0,$$

and

$$P(\{\mathbf{1}_A\mathcal{E}_g[\alpha\xi + (1-\alpha)\eta | \mathcal{F}_t] - \mathbf{1}_A(\alpha\mathcal{E}_g[\xi | \mathcal{F}_t] + (1-\alpha)\mathcal{E}_g[\eta | \mathcal{F}_t]) > 0\}) > 0.$$

Since $g(t, 0) \equiv 0$ and $A \in \mathcal{F}_t$, it is obvious that

$$(3.13) \qquad \mathcal{E}_g[\mathbf{1}_A X | \mathcal{F}_t] = \mathbf{1}_A\mathcal{E}_g[X | \mathcal{F}_t] \qquad \forall X \in L^2(\mathcal{F}_T).$$

Since $g$ is independent of $y$ and $A \in \mathcal{F}_t$, by Definition 1.2, Lemma 1.1(iv), equality (3.13) and Lemma 1.1(iii) we infer

$$\mathcal{E}_g[\mathbf{1}_A(\alpha\xi + (1-\alpha)\eta) - \mathbf{1}_A(\alpha\mathcal{E}_g[\xi | \mathcal{F}_t] + (1-\alpha)\mathcal{E}_g[\eta | \mathcal{F}_t])]$$

$$= \mathcal{E}_g\{\mathcal{E}_g[\mathbf{1}_A(\alpha\xi + (1-\alpha)\eta)$$

$$- \mathbf{1}_A(\alpha\mathcal{E}_g[\xi | \mathcal{F}_t] + (1-\alpha)\mathcal{E}_g[\eta | \mathcal{F}_t]) | \mathcal{F}_t]\}$$

$$(3.14) \qquad = \mathcal{E}_g\{\mathcal{E}_g[\mathbf{1}_A(\alpha\xi + (1-\alpha)\eta) | \mathcal{F}_t]$$

$$- \mathbf{1}_A(\alpha\mathcal{E}_g[\xi | \mathcal{F}_t] + (1-\alpha)\mathcal{E}_g[\eta | \mathcal{F}_t])\}$$

$$= \mathcal{E}_g\{\mathbf{1}_A\mathcal{E}_g[(\alpha\xi + (1-\alpha)\eta) | \mathcal{F}_t]$$

$$- \mathbf{1}_A(\alpha\mathcal{E}_g[\xi | \mathcal{F}_t] + (1-\alpha)\mathcal{E}_g[\eta | \mathcal{F}_t])\} > 0.$$

On the other hand, since $\mathcal{E}_g$ is convex and $g$ is independent of $y$, in view of Definition 1.2, Lemma 1.1(iv) and equality (3.13), we deduce that

$$\mathcal{E}_g[\mathbf{1}_A(\alpha\xi + (1-\alpha)\eta) - \mathbf{1}_A(\alpha\mathcal{E}_g(\xi|\mathcal{F}_t) + (1-\alpha)\mathcal{E}_g(\eta|\mathcal{F}_t))]$$
$$= \mathcal{E}_g[\alpha(\mathbf{1}_A\xi - \mathbf{1}_A\mathcal{E}_g(\xi|\mathcal{F}_t)) + (1-\alpha)(\mathbf{1}_A\eta - \mathbf{1}_A\mathcal{E}_g(\eta|\mathcal{F}_t))]$$
(3.15)
$$\leq \alpha\mathcal{E}_g[\mathbf{1}_A\xi - \mathbf{1}_A\mathcal{E}_g(\xi|\mathcal{F}_t)] + (1-\alpha)\mathcal{E}_g[\mathbf{1}_A\eta - \mathbf{1}_A\mathcal{E}_g(\eta|\mathcal{F}_t)]$$
$$= \alpha\mathcal{E}_g\{\mathcal{E}_g[\mathbf{1}_A\xi - \mathbf{1}_A\mathcal{E}_g(\xi|\mathcal{F}_t)|\mathcal{F}_t]\}$$
$$\quad + (1-\alpha)\mathcal{E}_g\{\mathcal{E}_g[\mathbf{1}_A\eta - \mathbf{1}_A\mathcal{E}_g(\eta|\mathcal{F}_t)|\mathcal{F}_t]\}$$
$$= \alpha\mathcal{E}_g\{\mathcal{E}_g[\mathbf{1}_A\xi|\mathcal{F}_t] - \mathbf{1}_A\mathcal{E}_g(\xi|\mathcal{F}_t)\}$$
$$\quad + (1-\alpha)\mathcal{E}_g\{\mathcal{E}_g[\mathbf{1}_A\eta|\mathcal{F}_t] - \mathbf{1}_A\mathcal{E}_g(\eta|\mathcal{F}_t)\}$$
$$= \alpha\mathcal{E}_g\{\mathbf{1}_A\mathcal{E}_g[\xi|\mathcal{F}_t] - \mathbf{1}_A\mathcal{E}_g(\xi|\mathcal{F}_t)\}$$
$$\quad + (1-\alpha)\mathcal{E}_g\{\mathbf{1}_A\mathcal{E}_g[\eta|\mathcal{F}_t] - \mathbf{1}_A\mathcal{E}_g(\eta|\mathcal{F}_t)\}$$
$$= 0 + 0 = 0.$$

Clearly (3.15) is a contradiction to (3.14). Therefore $P(A) = 0$. Thus (3.12) does hold.

For any $z_1, z_2 \in \mathbf{R}^d, \alpha \in [0,1]$, if $t \in S^{\alpha z_1 + (1-\alpha)z_2}(g) \cap S^{z_1}(g) \cap S^{z_2}(g)$, by (3.12) we deduce that

$$P\text{-a.s.,} \qquad g(t, \alpha z_1 + (1-\alpha)z_2) \leq \alpha g(t, z_1) + (1-\alpha)g(t, z_2).$$

For any $z_1, z_2 \in \mathbf{R}^d, \alpha \in [0,1]$, by the Representation Lemma we know that

$$\lambda([0,T] \setminus (S^{\alpha z_1 + (1-\alpha)z_2}(g) \cap S^{z_1}(g) \cap S^{z_2}(g))) = 0.$$

Thus for any $z_1, z_2 \in \mathbf{R}^d, \alpha \in [0,1]$, we have

(3.16)   $dP \times dt$-a.s.,    $g(t, \alpha z_1 + (1-\alpha)z_2) \leq \alpha g(t, z_1) + (1-\alpha)g(t, z_2).$

Thus (iii) does hold.   $\square$

Analogously to the argument of convexity for $g$-expectations, we have:

THEOREM 3.3 (Subadditivity for $g$-expectation). *Let* (A1) *and* (A3) *hold for* $g$. *Then the following three statements are equivalent:*

(i) $\mathcal{E}_g[\cdot]$ *is subadditive.*                              (Subadditivity.)

(ii) *For any* $t \in [0,T]$, $\mathcal{E}_g[\cdot|\mathcal{F}_t]$ *is subadditive, that is,* $\forall \xi, \eta \in L^2(\mathcal{F}_T)$,

$$\mathcal{E}_g[\xi + \eta|\mathcal{F}_t] \leq \mathcal{E}_g[\xi|\mathcal{F}_t] + \mathcal{E}_g[\eta|\mathcal{F}_t], \qquad P\text{-a.s.}$$

(iii) $g$ *is independent of* $y$ *and* $g$ *is subadditive with respect to* $z$, *that is, for any* $z_1, z_2 \in \mathbf{R}^d$,

$$g(t, z_1 + z_2) \leq g(t, z_2) + g(t, z_2), \qquad dP \times dt\text{-a.s.}$$

For the positive homogeneity of *g*-expectations, we have the following theorem.

THEOREM 3.4 (Positive homogeneity for *g*-expectation). *Let* (A1) *and* (A3) *hold for g. Then the following three statements are equivalent:*

(i) $\mathcal{E}_g[\cdot]$ *is positively homogeneous.*

(ii) *For any* $t \in [0, T]$, $\mathcal{E}_g[\cdot | \mathcal{F}_t]$ *is positively homogeneous, that is,* $\forall \xi \in L^2(\mathcal{F}_T)$, $\alpha \geq 0$,

$$\mathcal{E}_g[\alpha \xi | \mathcal{F}_t] = \alpha \mathcal{E}_g[\xi | \mathcal{F}_t], \qquad P\text{-}a.s.$$

(iii) *g is positively homogeneous with respect to* $(y, z)$, *that is, for any* $(y, z) \in \mathbf{R} \times \mathbf{R}^d$, $\alpha \geq 0$,

$$g(t, \alpha y, \alpha z) = \alpha g(t, y, z), \qquad dP \times dt\text{-}a.s.$$

PROOF. (iii) $\Rightarrow$ (ii) is just Proposition 9 of [21]. (ii) $\Rightarrow$ (i) is trivial. Now let us prove that (i) $\Rightarrow$ (iii). Suppose that (i) holds for $\mathcal{E}_g$. For any $\alpha > 0$, we define a new function

$$\tilde{g}^\alpha(t, y, z) := \alpha g\left(t, \frac{y}{\alpha}, \frac{z}{\alpha}\right) \qquad \forall (t, y, z) \in [0, T] \times \mathbf{R} \times \mathbf{R}^d.$$

It is clear that $\tilde{g}^\alpha$ satisfies (A1), (A2) and (A3).

For any $\xi \in L^2(\mathcal{F}_T)$, we deduce that

$$(Y_t(\tilde{g}^\alpha, T, \alpha\xi), Z_t(\tilde{g}^\alpha, T, \alpha\xi))_{t \in [0, T]} = \alpha(Y_t(g, T, \xi), Z_t(g, T, \xi))_{t \in [0, T]}.$$

Thus for any given $\alpha > 0$, we have

$$\mathcal{E}_{\tilde{g}^\alpha}[\alpha\xi] = \alpha\mathcal{E}_g[\xi] \qquad \forall \xi \in L^2(\mathcal{F}_T).$$

Combining (i) with the above equality we have

$$(3.17) \qquad \mathcal{E}_{\tilde{g}^\alpha}[\xi] = \mathcal{E}_g[\xi] \qquad \forall \xi \in L^2(\mathcal{F}_T).$$

Thus for any $\alpha > 0$, using the same argument as in (3.1)–(3.6) we conclude that $\tilde{g}^\alpha = g$, that is,

$$(3.18) \quad dP \times dt\text{-a.s.}, \forall (y, z) \in \mathbf{R} \times R^d, \qquad g(t, y, z) = \alpha g\left(t, \frac{y}{\alpha}, \frac{z}{\alpha}\right).$$

Hence (iii) follows from (A1) and (3.18). □

REMARK 3.2. Under (A1), (A3) and (A4), [21] proved that (ii) is equivalent to (iii) in Theorems 3.2–3.4; see Propositions 8–11 in [21].

**4. Risk measures via $g$-expectations.** Recently, many papers have been devoted to the problem of quantifying the risk of a financial position. Such a position, as in Artzner et al. [1] and Föllmer and Schied [10], will be described by the corresponding payoff profile, that is, by a real-valued function $X$ on some set $\Omega$ of possible scenarios, where $X(\omega)$ is the discounted net worth of the position at the end of the trading period if the scenario $\omega \in \Omega$ is realized. Coherent risk measures were introduced by Artzner et al. [1]; then, convex risk measures were introduced by Föllmer and Schied [9], and independently, by Frittelli and Rosazza Gianin [11]. Among others, we are especially interested in [1, 2, 7, 9, 10, 11, 19, 20, 21]. For the convenience of the reader, we recall some definitions of risk measures. The definition of convex measure of risk we use in this paper was given by Föllmer and Schied [10], slightly different from the one given by Frittelli and Rosazza Gianin [11].

DEFINITION 4.1 ([1, 10]).   Let $\mathcal{G}$ be the set of risks, that is, a set of real-valued functions on $\Omega$. A mapping $\rho : \mathcal{G} \to \mathbf{R}$ is called a monetary risk measure if $\rho$ satisfies the following conditions for all $X, Y, \in \mathcal{G}$:

(1) Monotonicity: If $X \leq Y$, then $\rho(Y) \geq \rho(X)$.

(2) Translation invariance: If $c \in \mathbf{R}$, then $\rho(X + c) = \rho(X) - c$.

DEFINITION 4.2 ([1, 10]).   A monetary risk measure $\rho$ is called a convex measure of risk if it satisfies

(3) Convexity: $\rho(\lambda X + (1 - \lambda)Y) \leq \lambda \rho(YX) + (1 - \lambda)\rho(Y)$, $\forall \lambda \in [0, 1]$.

A convex measure of risk $\rho$ is called a coherent measure of risk if it satisfies

(4) Positive homogeneity: If $\lambda \geq 0$, then $\rho(\lambda X) = \lambda \rho(X)$.

DEFINITION 4.3 ([21]).   Let $\mathcal{G}$ be the set of risks and $\mathcal{G} \subseteq L^0(\Omega, \mathcal{F}_T, P)$. A map system $(\rho_t)_{t \in [0,T]}$ is called a dynamic risk measure if it satisfies the following conditions for all $X, Y \in \mathcal{G}$ and all $t \in [0, T]$:

(1°) $\rho_t : \mathcal{G} \mapsto L^0(\Omega, \mathcal{F}_t, P)$.

(2°) $\rho_0$ is a static monetary risk measure.

(3°) $\rho_T(X) = -X$.

(4°) Dynamic monotonicity: If $X \leq Y$, then $\rho_t(Y) \geq \rho_t(X)$.

(5°) Dynamic translation invariance: If $c \in \mathbf{R}$, then $\rho_t(X + c) = \rho_t(X) - c$.

DEFINITION 4.4 ([21]).   A dynamic risk measure $(\rho_t)_{t \in [0,T]}$ is called a dynamic convex measure of risk if it satisfies

(6°) Dynamic convexity: For any $X, Y \in \mathcal{G}$, $\lambda \in [0, 1]$, $t \in [0, T]$, $P$-a.s.,

$$\rho_t(\lambda X + (1 - \lambda)Y) \leq \lambda \rho_t(X) + (1 - \lambda)\rho_t(Y).$$

A dynamic convex measure of risk $(\rho_t)_{t \in [0,T]}$ is called a dynamic coherent measure of risk if it satisfies

(7°) Dynamic positive homogeneity: For any $X \in \mathcal{G}$, $\lambda \geq 0$, $t \in [0, T]$,

$$P\text{-a.s.,} \qquad \rho_t(\lambda X) = \lambda \rho_t(X).$$

By Theorems 3.1–3.4 and Definitions 4.1–4.4, we can obtain the following Theorems 4.1 and 4.2 immediately.

THEOREM 4.1. *Let* (A1) *and* (A3) *hold for* $g$. *Let the set* $\mathcal{G}$ *of risks be* $L^2(\mathcal{F}_T)$. *Let* $\rho^g$ *and* $(\rho_t^g)_{t \in [0, T]}$ *be defined as in equality* (1.2). *Then the following statements are equivalent:*

 (i) $\rho^g$ *is a convex measure of risk.*
 (ii) $(\rho_t^g)_{t \in [0, T]}$ *is a dynamic convex measure of risk.*
 (iii) $\mathcal{E}_g$ *is convex.*
 (iv) $g$ *is independent of* $y$ *and is convex with respect to* $z$.

THEOREM 4.2. *Let* $g$, $\mathcal{G}$, $\rho^g$ *and* $(\rho_t^g)_{t \in [0, T]}$ *be as in Theorem* 4.1. *Then the following statements are equivalent:*

 (i) $\rho^g$ *is a coherent measure of risk.*
 (ii) $(\rho_t^g)_{t \in [0, T]}$ *is a dynamic coherent measure of risk.*
 (iii) $\mathcal{E}_g$ *is sublinear, that is,* $\mathcal{E}_g$ *is positively homogeneous and subadditive.*
 (iv) $g$ *is independent of* $y$ *and is sublinear with respect to* $z$.

REMARK 4.1. Rosazza Gianin [21] proved that (ii) is equivalent to (iv) in Theorems 4.1 and 4.2 under assumptions (A1), (A3) and (A4).

REMARK 4.2. Generally, verifying that a dynamic risk measure $(\rho_t)_{t \in [0, T]}$ is a dynamic convex (resp. coherent) measure may be much more difficult than verifying that the corresponding static risk measure $\rho_0$ is a static convex (resp. coherent) measure. But for risk measures $(\rho_t^g)_{t \in [0, T]}$, by Theorems 4.1 and 4.2 we know that if $\rho^g$ is a static convex (resp. coherent) measure of risk, then $(\rho_t^g)_{t \in [0, T]}$ must be a dynamic convex (resp. coherent) measure of risk.

**Acknowledgments.** This paper is an improved version of part of my Ph.D. thesis. I thank Professors S. Peng and Z. Chen for their help. I also thank Professor R. Adler and three anonymous referees for their careful reading of this paper and valuable suggestions.

Department of Mathematics
China University of Mining and Technology
Xuzhou 221116
P. R. China
E-mail: jianglong365@hotmail.com